\documentclass[12pt]{article}
\usepackage{amsmath,amssymb,amsthm}
\usepackage{mathtools}
\usepackage{graphics,epsfig}
\usepackage{hyperref}
\usepackage{natbib}
\usepackage{color}
\usepackage{graphicx}
\usepackage{caption}
\usepackage{subcaption}
\usepackage{float}
\usepackage{mathrsfs}
\usepackage{multirow}
\usepackage{comment}
\usepackage{setspace}
\usepackage{bbm}
\usepackage{bm}
\usepackage{amsfonts}
\usepackage{xcolor}
\usepackage{pslatex}
\usepackage{setspace}
\usepackage{bbm}
\def \b {\boldsymbol}

\def \E {\mathbb{E}}

\def \a {\alpha}
\def \be {\beta}


\begin{document}

  \title{\bfseries Optimal drug application on stochastic cancer growth: an approach through path integral control }

  \author{
Noelymar Farinacci\footnote{Corresponding author, {\small\texttt{nf2322@jagmail.southalabama.edu}}}\; \footnote{College of Nursing,  University of South Alabama, 5721 USA N Dr, Mobile, AL, 36688,
United States.}
}

\date{\today}
\maketitle

\begin{abstract}
We provide an overview of an optimal control problem within a stochastic model of tumor growth, which includes drug application. The model comprises two stochastic differential equations (SDE) representing the diffusion of nutrient and drug concentrations. To account for various uncertainties, stochastic terms are incorporated into the deterministic framework, capturing random disturbances. Control variables, informed by medical principles, are used to regulate drug and nutrient concentrations. In defining the optimal control problem, a stochastic cost function can be established, and a Feynman-type path integral control approach would lead to an optimal drug treatment.
\end{abstract}

\section*{Introduction:}
Cancer remains among the top causes of mortality globally, sparking extensive research across various scientific fields. Mathematicians in cancer research have approached tumor growth through mathematical modeling from multiple perspectives \citep{zhao2010parabolic}. Frequently, models assume that tumors expand with radial symmetry, a concept supported by early in vitro observations showing that initial solid tumor growth tends to be radially symmetric \citep{esmaili2024optimal}. Evidence also suggests that low glucose and oxygen levels in the central parts of tumor spheroids may foster sub-populations of quiescent, hypoxic, anoxic, and necrotic cells \citep{tao2009free,zhao2010parabolic}. Consequently, many tumor models differentiate live cells into proliferative and quiescent groups. It's essential to note that all systems encounter uncertainties due to environmental factors, experimental variability, and more. This has led researchers to explore stochastic tumor growth models to account for random fluctuations and unknowns. Additionally, given the significance of cancer treatment, numerous studies have focused on models addressing the effects of drug therapies on tumors \citep{esmaili2017optimal,esmaili2024optimal}.

The rise of personalized medicine in cancer treatment has transformed therapeutic strategies for patients with tumors that possess actionable mutations. For some, this shift has significantly extended lifespans, lowered toxicity, and enhanced quality of life. Unfortunately, however, this group remains limited; estimates from 2020 suggest that only around 5\% of patients benefit from targeted therapies \citep{wang2024threshold}. Moreover, while personalized treatments offer substantial benefits, they rarely lead to a complete cure, as tumors often evolve resistance through Darwinian mechanisms \citep{khan2024mp60,vikramdeo2024abstract}. In response, a new approach called \emph{evolutionary therapy} aims to leverage disease evolution dynamics to optimize drug choices and treatment schedules. Researchers are utilizing both mathematical and experimental models to explore practical questions of theoretical significance, such as how resistance to one drug influences sensitivity to others and whether heterogeneous tumor populations respond differently to treatment based on their state. Insights from these studies have already informed rational drug sequencing and cycling, benefiting bacterial infections as well as various cancers. In therapeutic scheduling, adaptive therapy—rooted in Evolutionary Game Theory (EGT) has shown encouraging results, not only in theory but also in a phase 2 trial for men with metastatic prostate cancer \citep{kaznatcheev2017cancer}. Experimentally, the principles of EGT have been validated in vivo, with ongoing development of in vitro assays and observations of interactive dynamics using these methods \citep{wang2024threshold}.

Cancer cells, like other populations of living organisms, consist of individual cells with unique behaviors and evolutionary backgrounds. Their interactions and life histories are filled with stochastic elements, including genetic variability, transitions in cell fate, differential drug responses, signaling differences, and micro-scale environmental fluctuations within the tumor. Many of these represent demographic stochasticity, where random variation can often be \emph{averaged out} if the population is large enough. This concept underpins the understanding of tumor heterogeneity by dividing cells into sub-populations \citep{pramanik2024dependence}. Such grouping is effective when mutation-selection dynamics are balanced so that only closely related genotypes with the same phenotype persist stably \citep{wang2024threshold}. These clusters, known as quasispecies, form distributions around a central genotype, with cells in each group displaying similar behavior despite minor genetic differences and random birth or death events. However, environmental stochasticity, which remains significant even in large populations as it encompasses random factors that affect entire groups simultaneously. In cancer, such fluctuations can result from factors unrelated to treatment, such as small but frequent changes in the physiology of the host. Moreover, individual cells within subpopulations will still respond differently to these broad disturbances. Therefore, implementation of \emph{environmental stochasticity} reflects an approach focused on modeling the averaged response of subpopulations to these overarching system-wide influences \citep{pramanik2024parametric}.

When modeling such fluctuations in continuous time, Stochastic Differential Equations (SDEs) often appear, with Stochastic Control Theory providing a means to optimize their behavior \citep{pramanik2021optimization,pramanik2021optimal,pramanik2023scoring}. This approach gives a mathematical perspective for making sequential decision making such as determining the appropriate drug dose over time while accounting for random shocks \citep{pramanik2020motivation}, like unpredictable shifts in the fitness of competing cancer cell subpopulations. Any static treatment approach will produce a random evolutionary path for the tumor and a random cumulative cost (which could represent total drug usage, time to recovery, or a blend of these factors) \citep{wang2024threshold}. Dynamic Programming (DP) \citep{pramanik2020optimization,pramanik2024optimization} aims to create equations for the cumulative cost of the best strategy, allowing that strategy to be expressed in feedback form, meaning treatment decisions (dose and duration) are updated according to the current state of the tumor rather than following a predetermined schedule. This approach can be implemented to various cancer models and treatment types, including strategies aimed at either stabilizing or eradicating the tumor. However, rather than choosing a strategy that minimizes the average expected cost, one can select the objective function that maximizes the probability of achieving a specific goal (e.g., meeting therapy objectives without surpassing a set cost limit). These risk-aware, or \emph{threshold-aware}, policies are designed to adapt the treatment plan based on the response of the tumor to previous drugs and the cost accumulated up to that point. 

 \cite{mckean1966class} dynamics is highly advantageous in stochastic control under random perturbations, particularly for systems where interactions are shaped by both individual behaviors and population-level influences \citep{carmona2013control,carmona2015forward}. First, McKean–Vlasov dynamics enables mean-field interactions \citep{lasry2007mean}, meaning that the behavior of a single cancer cell depends not only on its own state but also on the state distribution of all other cells \citep{polansky2021motif,pramanik2024estimation}. Second, in large systems, modeling each cancer cell individually can be computationally intensive. McKean–Vlasov dynamics reduces this complexity by focusing on the average effect of the population, making large-scale stochastic control problems more computationally feasible \citep{hua2019assessing}. Third, mean-field terms help smooth fluctuations due to random noise, fostering more stable and reliable treatment policies a major benefit for systems where high noise sensitivity could otherwise result in erratic responses \citep{pramanik2022lock,pramanik2022stochastic,pramanik2023path1,pramanik2023path}. Fourth, in multi-host environments where hosts are influenced by their surroundings and each other, McKean–Vlasov dynamics captures these interdependent behaviors effectively. It allows for both centralized and decentralized strategies to efficiently manage large populations of interacting hosts \citep{pramanik2024semicooperation}. Fifth, McKean–Vlasov dynamics is beneficial for risk-sensitive or mean-field treatment strategies, as it enables the management of not only individual outcomes but also their distributions. This ability to fine-tune policies to control both mean and variance is particularly useful when controlling risk (variance) is essential \citep{pramanik2016tail,pramanik2021effects}. Lastly, in systems exhibiting herding effects, McKean–Vlasov dynamics incorporates self-reinforcing behaviors, where hosts’ actions influence one another \citep{pramanik2024motivation}. This approach is valuable in settings requiring behavior synchronization or collective dynamics for the desired system outcome.

\section*{Formulation of a stochastic control:}
Let us assume a fixed finite time $t>0$ where the cancer treatment takes place. For a measurable space $(\Omega,\mathcal F)$, let $\mathcal P(\Omega,\mathcal F)$ be the probability measure on $(\Omega,\mathcal F)$ such that $\Omega$ is the sample space of the state of cancer, and $\mathcal F$ is the Borel sigma algebra $\mathcal B(\Omega)$ \citep{carmona2016mean}. Let $\mathfrak C^n=C\left([0,t]\big|\mathbb R^n\right)$ be the continuous functions $[0,t]\mapsto \mathbb R^n$. Define an evaluation process $\zeta_s$ on $\mathfrak C^n$ by $\zeta_s(x)=x(s)$ and the truncated supremum norms $||.||_s$ on $\mathfrak C^n$ by
\[
||x||_s:=\sup_{\nu\in[0,s]}|x_\nu|, \ \ \ \forall\ s\in[0,t],
\]
where $x$ is a value a value of the cancer state $X\subset\mathcal X$. Unless otherwise specified, $\mathfrak C^n$ is endowed with the norm $||.||_t$. Assume $\bm B(s)$ be an $n$-dimensional Brownian motion under $\mathfrak C^n$. For $\theta\in\mathcal P(\mathfrak C^n)$, let the image of $\theta$ under $\zeta_s$ is represented by $\theta(s)\in\mathcal P(\mathbb R^n)$. For $\rho\geq 0$ and a separable metric space $(\kappa,\gamma)$, assume $\mathcal P^\rho(\kappa)$ is the set of probability laws $\theta\in\mathcal P(\kappa)$ with $\int_\kappa\gamma^\rho(x,x_0)\theta(dx)<\infty$ for some $x_0\in \kappa$. For $\rho\geq 1$ and $\theta,\theta'\in\mathcal P(\kappa)$, consider $\gamma_{\kappa,\gamma}$ as the $\rho$-Wasserstein distance, defined as 
\begin{equation}\label{0}
 \gamma_{\kappa,\gamma}(\theta,\theta'):=\inf\left\{\left[\int_{\kappa\times\kappa}\vartheta(dx,dy)\gamma^\rho(x,y)\right]^{1/\rho} : \text{$\vartheta\in\mathcal P(\kappa\times\kappa)$ has marginals $\theta,\theta'$}\right\}.
\end{equation}
The probability space $\mathcal P^\rho(\kappa)$ equipped with the measure $\gamma_{\kappa,\gamma}$, and $\mathcal P(\kappa)$ has topology of weak convergence \citep{carmona2016mean,pramanik2024optimization}. Both of these probability spaces are defined on Borel $\sigma$-algebras, which coincide with the $\sigma$-algebra generated by the mappings
\begin{align*}
 &\mathcal P^\rho(\kappa)\ni\theta\mapsto \theta(\tilde\kappa),\\
& \mathcal P(\kappa)\ni\theta\mapsto \theta(\tilde\kappa),
\end{align*}
where $\tilde \kappa$ is an arbitrary subset of $\kappa$. For a fixed adapted therapy $\{\theta(s)\}_{s\in[0,t]}$ with values in $\mathcal P(\mathbb R^n)$ of probability measure from $\mathbb R^n$, our objective to solve
\begin{equation}\label{1}
 \inf_u\mathcal J (x_0,u):=\inf_u\ \E\left\{\int_0^th[s, X(s),u(s),\theta(s)]ds+\tilde h[X(t),\theta(t)]\right\},   
\end{equation}
subject to a McKean-Vlasov dynamics represented by the SDE
\begin{equation}\label{2}
 dX(s)=\mu[s, X(s),u(s),\theta(s)]ds+\sigma[s, X(s),u(s),\theta(s)]d\bm B(s),  
\end{equation}
where for two exponents $\rho',\rho\geq 1$, and the adaptive treatment space $\mathcal U$ taking values from $\mathbb R^n$ consider the following mappings
\begin{align}\label{3}
 & (\mu, h):[0,t]\times\mathbb R^n\times\mathcal U\times\mathcal P^\rho(\mathbb R^n) \mapsto\mathbb R^n\times\mathbb R,\notag\\
 &\bm\sigma:[0,t]\times\mathbb R^n\times\mathcal U\times\mathcal P^\rho(\mathbb R^n) \mapsto\mathbb R^n\times\mathbb R,\notag\\
 &\tilde h:\mathbb R^n\times\mathcal P^\rho(\mathbb R^n)\mapsto\mathbb R.
\end{align}
In Equation \eqref{2} $\mu$, $\bm\sigma$ are the drift and the diffusion components, respectively, $u\in\mathcal U$ is an adaptive treatment, $h$ is the treatment cost function with its terminal value $\tilde h$ which takes the values from two sets, failure and success of treatments. $\tilde h$ has the following construction,
\[
\tilde h[X(t),\theta(t)]=\begin{cases}
    \infty, & \text{if $x(t)\in\mathcal T$},\\
    0, & \text{if $x(t)\in\mathcal T'$},
\end{cases}
\]
such that if $x(t)\in \mathcal T$, the therapy fails, and if $x(t)\in \mathcal T'$, the therapy succeeds. Now we are going to provide an explicit structure of the drift and the diffusion components of the SDE \eqref{2}. Following \cite{wang2024threshold} one can examine three categories of cancer cells. \emph{Glycolytic cells} (GLY) operate without oxygen and generate lactic acid, which can harm the neighboring healthy tissue. The other two cell types are \emph{aerobic}, with access to improved vasculature, which is encouraged by the production of the VEGF signaling protein. Cells that actively produce VEGF (VOP) allocate part of their resources toward developing vasculature, whereas the remaining aerobic cells do not contribute to this process and can be considered free-riders or defectors in game theory terms. Let $\{m_g(s),m_d(s),m_v(s)\}'$ be a vector of continuous time-dependent subpopulations of these three types of cancer with dynamics $dm_j/ds=\pi_im_i$, for all $i\in\{g,d,v\}$ so that $\pi_i$ is the fitness of $i^{th}$ type of cancer. These $\pi_i$'s are expressed from the competition between different populations of cancer cells. In this model of cellular competition within tumors, the dynamics resemble a \emph{public goods} or \emph{club goods} scenario. Here, VEGF functions as a \emph{club goods} since its benefits are limited to VOP and DEF cells, whereas the acid produced by GLY acts as a \emph{public goods} with its harmful effects on healthy tissue assumed to provide an advantage to all cancer cells. The foundational model in \citep{kaznatcheev2017cancer} presumes that each cell engages with $n$ neighboring cells, with its benefits determined by its own cell type and the distribution of other cell types nearby. Given that all cells are sampled uniformly at random from a large, well-mixed population, it is possible to calculate all fitness values, denoted as $\pi_i$, as the expected payoffs within this interaction involving $(n+1)$ cells \citep{wang2024threshold}. Since, the total cell population is extremely large (i.e., $n \rightarrow \infty$), we can approximate their interactions as resembling a mean-field game \citep{lasry2007mean, carmona2016mean}. Furthermore, the expected payoffs $\pi_i$'s depend on current subpopulation fractions
\begin{equation*}
  b_g=\frac{m_g}{m_g+m_d+m_v},\ b_d=\frac{m_d}{m_g+m_d+m_v},\ \text{and}\ b_v=\frac{m_v}{m_g+m_d+m_v}.
\end{equation*}
Rewriting everything in terms of the proportion of glycolytic cells in the tumor $X_2(s)=b_g(s)$, and the proportion of VOP among aerobic cells $X_1(s)=b_v(s)/\left[b_v(s)+b_d(s)\right]$. Therefore, the explicit structure of of the McKean-Vlasov SDE in Equation \eqref{2} is 
\begin{align}\label{4}
  X&=\begin{bmatrix}
      X_1(s)\\X_2(s)
  \end{bmatrix}
  =\begin{bmatrix}
      \frac{b_v(s)}{b_v(s)+b_d(s)}\\b_g(s)
  \end{bmatrix}, \ 
  \bm B(s)=\begin{bmatrix}
      B_g(s)\\ B_d(s)\\B_v(s)
  \end{bmatrix},\notag\\
 \mu&= \begin{bmatrix}
     X_1(s)[1-X_1(s)]\left\{\left[\frac{\be_v}{n+1}\left(\sum_{j=0}^dQ^j\right)-c\right]+\theta(s)+(1-X_1(s))\sigma_2^2-X_1(s)\sigma_3^2\right\}\\[5ex]
     \begin{multlined}
     X_2(s)[1-X_2(s)]\biggr\{\left[\frac{\be_\a}{n+1}-(\be_v-c)X_1(s)-u(s)\right]+\theta(s)\\[-4ex]
     -\left[\sigma_1^2X_2(s)-\sigma_2^2(1-X_2(s))(1-X_1(s))^2-\sigma_3^2(1-X_2(s))X_1(s)^2\right]\biggr\}
     \end{multlined}
 \end{bmatrix},\notag\\
 \bm\sigma&=\begin{bmatrix}
     \sigma_{11}(s) & \sigma_{12}(s) & \sigma_{13}(s) \\
     \sigma_{21}(s) & \sigma_{22}(s) &\sigma_{23}(s)
 \end{bmatrix}
\end{align}
where 
\begin{align*}
    \sigma_{11}(s)&=0\\
    \sigma_{12}(s)&=-\sigma_d X_1(s)(1-X_1(s))+\theta(s)\\
    \sigma_{13}(s)&=\sigma_v X_1(s)(1-X_1(s))+\theta(s)\\
    \sigma_{21}(s)&=\sigma_g X_2(s)(1-X_2(s))+\theta(s)\\
    \sigma_{22}(s)&=\sigma_g X_2(s)(1-X_2(s))(1-X_1(s))+\theta(s)\\
    \sigma_{23}(s)&=\sigma_g X_2(s)(1-X_2(s))X_1(s)
\end{align*}
In the expressions in the system \eqref{4}, $\bm B(s)$ is a standard three-dimensional Brownian motion for GLY, DEF and VOP cells, with volatilities $\sigma_g,\sigma_v,\sigma_d$, respectively,$\be_\a$ is the benefit per acidification, $\be_v$ is the benefit from oxygen per unit of vascularization, $c$ is the cost per production of VEGF. Our objective is to minimize the treatment cost functional
\begin{equation}\label{5}
    \mathcal J\left[X_1(0),X_2(0),u\right]=\E\left\{\int_0^tu(s)ds+et+\tilde h(X_1(t),X_2(t),\theta(t)\right\},
\end{equation}
where $e>0$ represents the relative importance of stabilization of cancer growth and failure of the treatment. The stochastic Lagrangian \citep{ewald2024adaptation} of this system is
\begin{multline}\label{6}
    \mathcal L[s, X(s),u(s),\theta(s),\lambda(s)]=\E\biggr\{\int_0^t h[s, X(s),u(s),\theta(s)]ds+\tilde h[X(t),\theta(t)]\\
    +\int_0^t\left[X(s)-x_0-\int_0^s\left[\mu[\nu,X(\nu),u(\nu),\theta(\nu)]d\nu-\bm\sigma[\nu,X(\nu),u(\nu),\theta(\nu)]d\bm B(\nu)\right]\right]d\lambda(s)\biggr\},
\end{multline}
where $\lambda(s)$ is the Lagrangian multiplier. To deal with Equation \eqref{6}, a specialized variant of the stochastic Pontryagin principle such as the Feynman-type path integral technique introduced by \cite{pramanik2020optimization} and further developed by \cite{pramanik2024optimization} can be applied to derive an analytical solution for this system \citep{pramanik2024semicooperation}. When dealing with high-dimensional state variables and nonlinear system dynamics, as encountered in Merton-Garman-Hamiltonian SDEs, numerically constructing a Hamilton-Jacobi-Bellman (HJB) equation becomes highly complex \citep{pramanik2023optimal}. The Feynman-type path integral method effectively addresses this issue of dimensionality by providing a localized analytical solution \citep{pramanik2021consensus,pramanik2024bayes}. To apply this approach, a stochastic Lagrangian (i.e., Equation \eqref{6}) is first formulated for each continuous point within the interval $s \in [0, t]$, where $t > 0$. Then, this time interval is split into $k$ equal subintervals, and a Riemann measure associated with the state variable is defined for each subinterval. After constructing a Euclidean action function, a Schr\"dinger-type equation is derived via Wick rotation \citep{djete2022mckean}. Imposing first-order conditions on both state and control variables allows us to find a solution of this system. This approach has promising applications in cancer research \citep{dasgupta2023frequent,hertweck2023clinicopathological,kakkat2023cardiovascular,khan2023myb,vikramdeo2023profiling,khan2024mp60}.

\section*{Concluding remarks:}
It is now well understood that cancers adapt and evolve while one person is under therapy, and this understanding is gradually being factored into treatment planning. In some instances, adapting treatment can be as straightforward as switching from one targeted therapy to another. However, instances where tumors consist of a diverse mix of interacting cell types, such adjustments are often impractical. Therefore, an ecological approach called \emph{adaptive therapy} is gaining interest. Until recently, clinical trials and theoretical models of adaptive therapy have relied on initial assumptions about the interactions within tumors and how they affect tumor composition over time \citep{wang2024threshold}. However, a growing number of studies, conducted both in lab settings and live models, are beginning to establish methods to measure more accurately these interactions. Since these tools advance, the next hurdle will be to understand these dynamics directly in patients and leverage them to enhance personalized treatment options.

 In this paper we suggest an alternative mean-field approach driven by a McKean-Vlasov dynamics. Mean field games, while present in various forms for some time, particularly in economics, have their theoretical roots in the groundbreaking contributions of \cite{lasry2006jeux1,lasry2006jeux2,lasry2007mean}, as well as \cite{huang2006large}. The central idea is to model the behavior of a large group of optimizing individuals who interact with each other through mean-field dynamics, under certain energetic or economic constraints. This approach simplifies the study of consensus by focusing on a control problem for a single representative individual, who interacts with the environment formed by the collective behavior of others. Essentially, when consensus is achieved, the system's inherent symmetries \citep{pramanik2024measuring} drive individuals to conform to a type of law of large numbers, resulting in a \emph{propagation of chaos} effect as the population size increases \citep{carmona2016mean}. Inspired by the theory of propagation of chaos, it can be anticipated that simplified equations might apply as the number of cancer cells becomes extremely large, making it feasible to work with solutions that are less complex than seeking Nash equilibria in large-scale stochastic differential games as previously discussed. This insight forms the basis for defining the mean field game problem. A mean field game is fundamentally an asymptotic framework for the game, where each cell's influence on the overall empirical measure is minimal. This suggests that, in the limit, optimization problems may become decoupled and uniform. In other words, the equilibrium problem at this limit reduces to a standard optimization problem for a single representative cancer cell, which interacts or competes solely with the environment shaped by the asymptotic behavior of the marginal empirical measures, $\{\theta(s)\}_{s\in[0,t]}$, linked to an adaptive treatment $\{u(s)\}_{s\in[0,t]}\in\mathcal U$ which is exchangeable. In the absence of common noise, the classical law of large numbers indicates that the limiting environment would be a deterministic sequence of probability measures, $\{\theta(s)\}_{s\in[0,t]}$, representing the equilibrium distribution of the population. However, when common noise is introduced, its effects do not cancel out, and thus the limiting environment becomes a stochastic sequence, $\{\theta(s)\}_{s\in[0,t]}$, reflecting the conditional distribution of the population in equilibrium given the common noise. This inclusion of common noise introduces new possibilities in cancer research by providing a dynamic perspective on the system's behavior.

\bibliographystyle{apalike}
\bibliography{bib}
\end{document}